\documentclass[11pt,a4paper]{article}
\newcommand{\bm}[1]{\mbox{\boldmath $#1$}}
\usepackage{epsfig}
\usepackage {amsfonts}
\usepackage {amsmath}
\usepackage{latexsym}

\begin{document}
\title{Construction of the similarity matrix for the spectral clustering method: numerical experiments}
\author{Paola Favati\thanks{IIT - CNR Via G. Moruzzi 1, 56124 Pisa, Italy ({\tt paola.favati@iit.cnr.it}).}
        \and Grazia Lotti\thanks{Dipartimento di Scienze Matematiche, Fisiche e Informatiche, University of Parma,
Parco Area delle Scienze 53A, 43124 Parma, Italy ({\tt grazia.lotti@unipr.it}).}\and
 Ornella Menchi\thanks{Dipartimento di
Informatica, University of Pisa, Largo Pontecorvo 3, 56127 Pisa,
Italy ({\tt menchi@di.unipi.it}).}  \and
 Francesco Romani\thanks{Dipartimento di
Informatica, University of Pisa, Largo Pontecorvo 3, 56127 Pisa,
Italy ({\tt romani@di.unipi.it}).} }

\date { }

\maketitle
\begin{abstract}
Spectral clustering is a powerful method for finding structure in a
dataset through the eigenvectors of a similarity matrix. It often
outperforms traditional clustering algorithms such as $k$-means when
the structure of the individual clusters is highly non-convex. Its
accuracy depends on how the similarity  between pairs of data points
is defined. Two important items contribute to the construction of the similarity matrix:  the sparsity of the underlying weighted graph, which depends mainly on the distances among data points, and the  similarity function.
When a Gaussian similarity function is used, the choice
of the scale parameter $\sigma$ can be critical.
In this paper we examine both items, the sparsity and the selection of suitable $\sigma$'s, based either directly on the  graph associated to the dataset or  on the minimal spanning tree (MST) of the graph.
An extensive  numerical experimentation
on  artificial and real-world datasets has been carried out to compare the performances of the methods.

\end{abstract}


\section{Introduction}\label{intro}
Clustering, a key step for many data mining problems, can be applied
to a variety of different kinds of documents as long as a distance
measure can be assigned to define the similarity among the data
objects. A clustering technique classifies the data into groups,
called {\it clusters}, in such a way that the objects belonging to a
same cluster are more similar to each other than to the objects belonging to different clusters. Clearly, a specific formulation of the clustering
problem depends considerably on the metric that is assigned.

We assume here that the data objects are the $n$ distinct points ${\bm
x}_1,\ldots,{\bm x}_n$ of a subset ${\cal X}$ of ${\bf
R}^m$. It is standard to measure the distances in ${\bf R}^m$ by using
norms. Both the 1-norm, known as Manhattan, and the 2-norm, known as
Euclidean, are suitable, while the $\infty$-norm is not (see
\cite{agga}, Sec. 3.2.1.4). We use
the Euclidean distance $d({\bm x}_i,{\bm x}_j)=\|{\bm x}_i-{\bm
x}_j\|_2$.
The given data objects are represented in a natural way by an {\it undirected  weighted  complete graph} $\mit\Gamma$, having  $n$  nodes and the distances $d({\bm x}_i,{\bm x}_j)$ as weights on the edges.
We assume also that the required number $k$ of clusters  is a
priori fixed, with $1<k\ll n$. Thus, a clustering ${{\mit
\Pi}=\{\pi_1,\ldots,\pi_k\}}$ corresponds to a partitioning of the
indices $\{1,\ldots,n\}$.

To solve the
clustering problem many different methods which work directly on the space ${\cal X}$ have been devised, but they are not entirely  satisfying. For example, the widely used
Euclidean-based $k$-means is very sensitive to the initialization and, applied directly to the original data,
cannot deal with clusters that are nonlinearly separable in ${\cal
X}$, has difficulties in taking into account other features like
the cardinality or the density of each cluster, and in some situations can be very slow \cite{vatt}. For this reason another approach can be followed using a graph-based procedure, which performs
the clustering by constructing the so called {\it similarity graph} $G$, having the same nodes of $\mit\Gamma$ and  possibly a subset of the edges of  $\mit\Gamma$. The weights on the edges of $G$ are defined by using a {\it similarity function}  to model the local neighborhood relationships between the data points, i.e. a higher similarity is assigned to pairs of closely related objects than to objects which are only weakly related. The $n\times n$ nonnegative symmetric adjacency matrix $W$ associated to $G$ is the so-called {\it similarity matrix}.
By using a spectral clustering algorithm the subspace ${\cal X}$,
where the original data live, is mapped to a new subspace ${\cal Y}$. The
transformed representation gives better results than the original
one, because the corresponding procedures are more flexible, can adjust to varying local densities and discover
arbitrarily shaped clusters (\cite{agga},
Sec. 6.7.1).
Then the clustering in the subspace ${\cal Y}$ is computed and mapped back to ${\cal X}$. Spectral clustering algorithms have become very popular (see \cite{na} for a thorough survey).
In the experiments we choose the algorithm described in \cite{ng}.

The similarity graph $G$ changes depending whether only distances among points are accounted for or density considerations are applied (see the literature review in \cite{inka2} for a list of the different approaches). In general the procedures which exploit specifically the density are more expensive.

A similarity function widely used in literature to construct the similarity matrix is the Gaussian one. It depends on a scale parameter $\sigma$ whose  choice can be crucial.
A possible way to
select $\sigma$ suggests to run the spectral algorithm repeatedly for
different values of  $\sigma$ and select the one that provides
the best clustering according to some criterium. Besides the obvious increase in the computational cost, the drawback of this procedure is the difficulty of choosing a reliable quality measure.

In this paper we take into consideration procedures for setting suitable values to $\sigma$ without performing multiple runs. These procedures use techniques based either directly on a subgraph $\mit\Delta$ of  $\mit\Gamma$ associated to the dataset or on the minimal spanning tree (MST) of $\mit\Delta$.
We focus on aspects like the sparsity level of the similarity matrix and the computational cost of its construction.
The paper is
organized as follows: after a very concise review of the
spectral clustering method, of the performance indices and of the minimum spanning tree in Section \ref{pre},
in Section \ref{graph} the different types of sparsity of the similarity matrices, the role of a scale parameter in their
construction and the different resulting methods we take into consideration are outlined.
Finally, in Section \ref{expe} the results
of the experimentation on both artificial and real-world datasets are presented and discussed.

\section{Preliminaries}\label{pre}

The idea behind a graph-based procedure is to convert the
representation of $n$ objects given in the {\it original} space
${\cal X}$ into a transformed representation in a new
$k$-dimensional {\it feature} space ${\cal Y}$, $k$ being the required number of clusters, in such a way that
the similarity structure of the data is preserved. To this aim, the
clustering problem is reformulated making use of the  similarity graph $G$.  The weight on the edge connecting the $i$th and the $j$th nodes
depends mainly on their
distance $d({\bm x}_i,{\bm x}_j)$ and is large for close points, small for
far-away points. Other aspects, like for example the nearby density, could be taken into consideration. For a correct application  of the spectral algorithm, the graph
$G$ should be undirected and connected \cite{lux}. Section \ref{graph} is dedicated to the description of the similarity graphs we consider.

\subsection{The normalized spectral clustering algorithm}\label{unn}
Once the  $n \times n$ similarity matrix $W$ associated to graph $G$ has been constructed, the spectral clustering can be performed. Spectral clustering algorithms are based on the eigen-decomposition of Laplacian matrices. A detailed description of many such algorithms can be found in \cite{na}.

The unnormalized  {\it Laplacian} matrix is $L=D-W$, where $D$ is the {\it degree} matrix of $W$, i.e. the
diagonal matrix whose principal entry is given by $d_{i,i}=\sum_{j=1}^n w_{i,j}$. Matrix $D$
can be seen as a local average of the similarity and is used as a
normalization factor for $W$. An exposition of the properties of matrix $L$ and of its normalized versions, on which spectral algorithms rely, can be found in the tutorial
\cite{lux}.
In our experimentation we use the algorithm proposed by \cite{ng} which exploits the normalized Laplacian
\[L_N=D^{-1/2}WD^{-1/2}.\]
This algorithm runs as follows:
\begin{enumerate}
\item construct matrix $L_N$,
\item compute the matrix $U=[{\bm u}_1$,
${\bm u}_2$,\ldots, ${\bm u}_k]$, where the ${\bm u}_i$ are the
mutually orthogonal eigenvectors of $L_N$ which correspond to the
largest $k$ eigenvalues (ordered downward).
\item perform a row normalization on $U$ to map the points ${\bm u}_i$ to a
unit hypersphere, by setting
 \[
  y_{i,j}=u_{i,j}\big/ \big(\sum_{r=1,k} u_{i,r}^{\,2}\big)^{1/2},
  \quad\hbox{for}\quad i=1,\ldots,n,\ j=1,\ldots,k.\]
For $i=1,\ldots,n$, points ${\bm y}_i$ of components $y_{i,j}$,  $j=1,\ldots,k$, are obtained.
The new feature space ${\cal Y}$, which is  a subspace of ${\bf R}^k$, is so identified.
\item compute the clustering of the points of ${\cal Y}$, which
is mapped back to the original
points of ${\cal X}$.
\end{enumerate}

The graph-based approach allows discovering
clusters of arbitrary shape. This property
is not specific of the particular clustering method used in
the final phase 4. Since the transformed representation enhances the
clustering properties of the data, even the simple  $k$-means can
be successfully applied. In our experimentation we apply the
clustering algorithm described in \cite{yu}, which is based on
a sequence of rotations and results to be fast and efficient.

Since the spectral clustering algorithm performs the dimensional
reduction through the eigenvectors of $L_N$, the
computation of the $k$ eigenvectors (using Lanczos
algorithm) can have substantial computational complexity for large
$n$ and $k$. The sparsity of matrix $W$ is then an important issue from this point of view.

\subsection{The performance index}\label{indi}
To measure the clustering validity one can rely on both internal or external
{\it indices} (see \cite{bern} and its references).
Internal indices generally evaluate the clustering in terms of compactness  and  separability, measured through Euclidean distances. When the data belong to $\mathbf{R}^2$ and $\mathbf{R}^3$ and the clusters are not convex or are badly separated or have different densities, these indices might not comply with what the human eye  would suggest.

External indices evaluate  how well the obtained clustering ${\mit
\Pi}_{obt}$ matches  the {\it target}, i.e. the assumed clustering ${\mit \Pi}_{tar}$ of the data.
The results obtained in the experimentation of Section \ref{expe} have been evaluated by means of the external {\it Normalized  Mutual Information index} (NMI) \cite{str1}. Three other external indices have also been used for comparison purpose, namely  the {\it purity} index \cite{man}, the {\it Rand} index \cite{ra} and the {\it Clustering Error} \cite{baj}. Since their grading on the set of all the methods that will be introduced in Section \ref{graph} was fully matching, in Section \ref{expe} we report only the clustering evaluations based on the NMI index.

Let ${{\mit \Pi}_{tar}=\{\widehat\pi_1,\ldots,\widehat\pi_k\}}$ be the  target  clustering, ${{\mit
\Pi}_{obt}=\{\pi_1,\ldots,\pi_{nc}\}}$ the computed clustering of $nc$ components and
define  $n_i=\# \;\pi_i$, $\widehat{n}_j=\#\; \widehat \pi_j $ and $n_{i,j}=\#\,(\pi_i\bigcap \widehat \pi_j)$, for $i=1,\ldots,nc$ and $j=1,\ldots,k$. The number $n_{i,j}$ counts how many points in the $i$th
cluster of ${\mit \Pi}_{obt}$  belong also to the $j$th cluster of
${\mit \Pi}_{tar}$. The NMI is given by
\[{\rm NMI}= \frac{2\,\sum\limits _{i=1}^{nc} \sum\limits_{j=1}^k n_{i,j} \log\big(\dfrac{n \ n_{i,j}}{n_i \ \widehat{n}_j}\big)}{ \sum\limits_{j=1}^k \widehat{n}_j \log \dfrac{\widehat{n}_j}{n} +\sum\limits_{i=1}^{nc} n_i \log \dfrac{n_i}{n} }.\]
In the literature there is also a version where the normalization is performed with the geometrical mean of the entropies instead of the arithmetical one. In the case $nc< k$, the algorithm is considered failed.

The larger the value NMI, the better the clustering performance according to the assumed target.
Perfectly matching clusterings have NMI=1, but lower indices can still be considered acceptable when the difference is due to few far-away points, possibly {\it outliers}. The treatment of outliers is not banal, as can be seen from the last three datasets in Figure \ref{fig1}, where it is not easy to decide whether the points of the two bridges should be considered outliers or belonging to the clusters of the other points (and which one) or even should form two separated clusters.
For this reason, no specific treatment  for the outliers detection is performed in the experiments.

\subsection{Minimum spanning tree}\label{mst}
Given a  connected  undirected weighted graph with $n$  nodes, its
 minimum spanning tree (MST)  is a subgraph
which has the same nodes, but connected by no more than $n-1$ edges
totalling the minimal weighting. Classical Prim's algorithm \cite{cor} computes
the MST with a computational cost  of order  $O(n^2)$ if the graph is dense and of order $O(e+n\log n)$, $e$ being the number of edges, if the graph is sparse, provided that  adjacency lists are used for the graph representation and a Fibonacci  heap is used as a working data structure \cite{cor}.

Prim's
algorithm starts from an arbitrary root node, constructs a tree
which grows until it spans all the nodes and returns the $n$-vector
${\bm p}$ such that $p_i$ is the parent of the $i$th node. MST is unique if there are no ties in the pairwise distances.
Since the graph is connected, also its MST is connected. From this property it follows that the graph obtained by removing  all edges longer than the longest edge of its MST is still connected \cite{lux}.

Clustering procedures which exploit the MST of a graph have been proposed \cite{gry,zho}. They act by removing from the graph the edges with the largest distances, until only $k$ disconnected components are obtained, avoiding to
drop edges which lead to subsets of a single element considered outliers. Anyway, MST based clustering algorithms are known to be highly instable in presence of noise and/or outliers \cite{sla}.

We are interested in MST not as a clustering procedure, but as a technique for extracting suitable values of the parameters that will be used in the construction of the similarity graphs, as described in the next section.

\section{The methods}\label{graph}
Sparsity is an important issue from the computational point of view: a sparse graph can be represented by adjacency lists which allow both a reduction of the required storage  and the implementation of more efficient matrix-vector product algorithms.
Hence, when $n$ is large, a sparse representation of the similarity among the data conveying sufficient information to get an acceptable clustering should be used.

The more immediate sparsity structures are obtained from {\it distance-based} graphs which take into account the pairwise distances $\delta_{i,j}=d({\bm x}_i,{\bm x}_j)$ of the $n$ points ${\bm x}_i$ in ${\bf R}^m$.
The sparse graph  obtained by dropping from the complete graph ${\mit\Gamma}$ some selected edges is denoted by ${\mit\Delta}$  and the similarity graph, which has the same nodes and edges of ${\mit\Delta}$ and weights defined through a similarity function, is denoted by $G$. The two cases will be considered:

\smallskip\noindent
$\bullet$ \quad the full case when ${\mit\Delta}={\mit\Gamma}$, i.e. no edge is dropped;

\smallskip\noindent
$\bullet$ \quad  the sparse case where  ${\mit\Delta}$ is a proper subgraph of ${\mit\Gamma}$ obtained according to a chosen sparsity model.

The  ratio $\nu/n^2$, where $\nu$ is the number of dropped edges from ${\mit\Gamma}$, measures the sparsity level. The greater $\nu$, the higher the sparsity.

\subsection{The sparsity model}\label{dbs}
In this paper we consider the distance-based sparsity obtained through techniques like the ones outlined in \cite{lux} where the sparsity level is imposed through a sparsity parameter.
The two following techniques are here considered.

\smallskip\noindent
(a)\quad $\epsilon$-{\it neighbor sparsity}, where $\epsilon$ is a selected threshold. In the sparse graph ${\mit \Delta}$ the nodes $i$ and $j$ are connected only if their pairwise distance $\delta_{i,j}$ is not greater than $\epsilon$.
As noted in \cite{lux}, if $\epsilon$ is set equal to the largest weight of the MST of ${\mit\Gamma}$, then graph $G$  results  connected, but this choice could give a too large  $\epsilon$ if some very far-away outliers are present.  For this reason smaller values should be chosen.  We denote these methods as $E$-methods.

\smallskip\noindent
(b)\quad $K$-{\it nearest neighbor sparsity}, where $K$ is the sparsity parameter. In the sparse directed
graph node $i$ is connected to node $j$ if ${\bm x}_j$ is among the $K$-nearest neighbors of ${\bm x}_i$.
Since the nearest neighborhood relationship is not symmetric, a symmetrization step must be provided to get an undirected graph. The symmetrization can be performed by connecting nodes $i$ and $j$ with an undirected  edge
\begin{description}
\item{(b1)}  if both points ${\bm x}_i$ and ${\bm x}_j$ are among the
$K$-nearest neighbors of the other one. We denote these methods as $M$-methods (the $M$ stands for mutual).
\item{(b2)} if either point
${\bm x}_i$ or point ${\bm x}_j$ is among the $K$-nearest
neighbors of the other one. We denote these methods as $N$-methods (the $N$ stands for non mutual).
\end{description}
Contrary to $\epsilon$-neighbor sparse graphs, $K$-nearest neighbor sparse graphs allow an  a-priori control of the achievable sparsity level. Moreover, $K$-nearest neighbor sparse graphs allow an easier detection of outliers. In fact, if before symmetrization no directed edge exists  from the node $j\ne i$ to the node $i$ for all $j$, the point ${\bm x}_i$ could be recognized as an outlier.

If the similarity graph obtained applying techniques (a) and (b) results to be unconnected, because too many connections have been lost, an aggregation step must be performed to guarantee sufficient connection.

For what concerns the $E$-methods, in the experimentation we set $\epsilon$ equal to the mean of the distances of each point from its $K$th closest neighbor.
In any case the chosen  integer $K$  depends on $n$ but is much smaller than $n$.
In \cite{lux} the value  $K_\ell=1+\lfloor\log_2 n\rfloor$ is suggested, on the basis of asymptotic connectivity considerations. In our experiments we consider also the value $K_s=1+\lfloor\sqrt n\rfloor$ which leads to a lower sparsity of the similarity graph without a significant increase of the computational cost for the similarity matrix construction (see Subsection \ref{cost}).

\subsection{The similarity function}\label{simfu}
Let ${\mit \Delta}_e$ denote the set of edges of the sparse graph ${\mit\Delta}$ obtained according one of the techniques described above. Obviously the edge $(i,i)\not\in {\mit \Delta}_e$ and $\delta_{i,j}\ne 0$ for $(i,j)\in {\mit \Delta}_e$.
A simple definition of the similarity function would be $w_{i,j}=1/\delta_{i,j}$ for $(i,j)\in {\mit \Delta}_e$, but this function is not satisfactory in presence of very close points.
Similarity functions considered in literature are:

\smallskip\noindent
the unit similarity  function
 \begin{equation}\label{uni}
w_{i,j}= 1\ \hbox{for}\ (i,j)\in {\mit \Delta}_e,
 \end{equation}
and the Gaussian similarity  function
 \begin{equation}\label{gau}
 w_{i,j}=
  \exp  \Big(-
   \dfrac{\delta_{i,j}^2}{2 \sigma^2}\Big)\ \hbox{for}\ (i,j)\in {\mit \Delta}_e,
   \end{equation}
where the {\it scale} parameter $\sigma$ controls the decay
rate of the distances.

The values $w_{i,j}$ are the weights of the similarity graph $G$.
The similarity matrix $W$ is the adjacency matrix of $G$, hence its entries are $w_{i,j}$ for $(i,j)\in {\mit \Delta}_e$ and zero otherwise.

\medskip
When the Gaussian similarity  function (\ref{gau}) is used, a correct choice of $\sigma$ may be critical for the efficiency of the spectral algorithm. A too small value of $\sigma$ would
give weights $w_{i,j}$ very close to 0, so that all the points would
appear equally far-away. On the contrary, a too large value of $\sigma$ would give weights $w_{i,j}$ very close to 1, so that all the
points would appear equally close. In both cases it would be
difficult to discriminate between close and distant points. Hence an intermediate value between the smallest and the largest $\delta_{i,j}$, with $i\ne j$, must be chosen.

\medskip
Instead of a {\it global} $\sigma$, different {\it local}
scales $\sigma_i$ can be used for the different objects ${\bm x}_i$.
Then (\ref{gau}) is replaced by
 \begin{equation}\label{gau1}
   w_{i,j}=
\exp  \Big(-
   \dfrac{\delta_{i,j}^2}{2 \sigma_i\sigma_j}\Big)\ \hbox{for}\ (i,j)\in {\mit \Delta}_e,
   \end{equation}
in order to
tune pairwise distances according to the local statistics of the
surrounding neighbors \cite{perona}. These local scales are suggested especially for
high-dimensional problems with large behavior variations, where the choice (\ref{gau1}) is trusted in better results than (\ref{gau}).

\medskip
A first simple idea for determining a reasonable value of $\sigma$ takes into account the decay rate of the exponential function $f(x)=\exp(-x^2/(2 \sigma^2))$ whose inflection point $x=\sigma$
discriminates between a first part of rapid decay and a second
flatter part. This suggests to look for gaps in the curve of
the distances $\delta_{i,j}$ sorted in descending order. If we can
individuate a large gap, we can exploit it as a suitable value for $\sigma$.

Another idea takes into account the histogram of the distances
$\delta_{i,j}$. If the data form clusters, the histogram is multi-modal and, if $\sigma$ is chosen around the first mode, the affinity values of the points which form a cluster
can be expected to be significantly larger than others. This
suggests to choose for $\sigma$ a value somewhat smaller than the
first mode of the histogram. In \cite{fisc} this
technique is said to work well for spherical-like clusters.

Although simple, these ideas are difficult to implement, because of the difficulty of detecting the right gap or mode. We have verified that their applicability is in effect
restricted to a small number of cases, and for this reason we have discarded them.

To
find a suitable value for the scale parameter, it is suggested in literature to run the spectral
clustering algorithm repeatedly for different values of $\sigma$ and
select the one that provides the best clustering according to a
chosen quality measure. For example, in \cite{ng} it is suggested to choose for $\sigma$ the value which gives the tightest clusters in ${\cal Y}$.
To implement this technique, a decreasing sequence could be used, starting with a value of $\sigma$ smaller than $\max_{i,j}\delta_{i,j}$. However, the drawback of this
procedure is the choice of the quality measure to use, which might be a non-monotone function of $\sigma$. In this case a reliable
determination of an acceptable $\sigma$ could not be obtained using
only a small number of tries. Moreover,
a test based on an internal performance index, like for example the DB index \cite{DB}, does not guarantee that the final clustering is consistent with what a human would set as a target.

\medskip
For a single run of the spectral clustering algorithm, a careful a-priori selection of $\sigma$ should be considered and we suggest in the next subsection to exploit either the  minimum spanning tree of $\mit \Delta$, denoted by MST$(\mit\Delta)$, or directly $\mit \Delta$. In the first case MST$(\mit\Delta)$ is assumed as a reliable representation of the data and $\sigma$ is set equal to its largest weight. In the second case, we suggest for  local $\sigma_i$'s the distance from ${\bm x}_i$ of its $\ell$th nearest neighbor, $\ell$ being a chosen index much smaller than $n$. In \cite{perona} the index $\ell=7$ is suggested, but this choice, independent from $n$, appears somewhat arbitrary. Our choice of  $\sigma_i$'s will be described in the next Section.
From local scales a single global $\sigma$ can be obtained by averaging the $\sigma_i$'s.

If the similarity  functions (\ref{gau}) or (\ref{gau1}) are used, the  graph $G$ might result to be sparse in practice also when no particular sparsity structure is explicitly imposed on ${\mit\Gamma}$, depending on the magnitude of the scale parameters, due to a possible fast decay of the exponential function. So, when  edges with negligible weights are dropped and the adjacency lists are used to represent  $G$, for all computational purposes this graph can be considered sparse.

\subsection{Definition of the methods}\label{def}
By the term ``method'' we intend the combination of a sparsity model and a similarity  function (either unit, global Gaussian or local Gaussian). The different  similarity matrices $W$ so constructed are then processed by the spectral algorithm.

\smallskip\noindent
$\bullet$\quad In the full case, no sparsity is imposed on ${\mit \Gamma}$, i.e.  ${\mit \Delta}={\mit \Gamma}$. In the following we refer to these methods as $F$-methods. The unit similarity function is not used.

\smallskip\noindent
--- \quad Method $F_1$ implements the similarity function (\ref{gau}). In \cite{lux} the largest weight of MST$({\mit \Gamma})$ is suggested as $\sigma$. In our experimentation this choice appeared often too large, matrix $W$ resulted nearly full and gave really poor clusterings. Hence if $\sigma$ results larger than the mean $\delta_{\rm mean}$ of all the pairwise distances ${\delta}_{i,j}$, we set $\sigma$ equal to $\delta_{\rm mean}$.

\smallskip\noindent
--- \quad Method $F_2$ implements the similarity  function (\ref{gau1}). The $i$th local scale $\sigma_i$ is set equal to the distance from ${\bm x}_i$ of its $K$th  nearest neighbor, with $K=K_\ell$. The sparsity level obtained with this value of $K$ is too low in most cases of the datasets considered in the experimentation, so the value $K=K_s$, which would give lower sparsity levels, has not been used.

\smallskip\noindent
--- \quad Method $F_3$ implements the similarity  function (\ref{gau}) with $\sigma$ equal to the mean of the previous $\sigma_i$'s.

\smallskip\noindent
$\bullet$\quad In the sparse case, let ${\mit\Delta}_e$ and MST$(\mit\Delta)_e$ be the edge sets of  $\mit\Delta$ and of MST$(\mit\Delta)$, respectively. We set
\begin{equation}\label{sigma}
\begin{array}{l}
t=\max\limits_{\textstyle{(i,j)\in \hbox{MST}{(\mit\Delta)}_{e}}}\ \delta_{i,j},\\[0.4cm]
s_i= \max\limits_{\textstyle{(i,j)\in {\mit\Delta}_e}}\ \delta_{i,j},\ \hbox{for}\ i=1,\ldots,n,\\[0.4cm]
s={\rm mean\ of\ the\ } s_i.
\end{array}\end{equation}
The name of a method indicates the sparsity model with a subscript which indicates the similarity function: specifically,  subscript 1 indicates similarity (\ref{uni}), subscript 2 indicates similarity (\ref{gau}) with $\sigma=t$, subscript 3 indicates similarity (\ref{gau1}) with $\sigma_i=s_i$, subscript 4 indicates similarity (\ref{gau}) with $\sigma=s$. For example $E_1$ denotes the $E$-method applied with the unit similarity function and $M_3$ denotes the $M$-method applied with local Gaussian function (\ref{gau1}).
Each method is applied with both $K=K_\ell$ and $K=K_s$.

\smallskip
On the whole, the  methods considered for the experiments are $q=27$.

\subsection{Computational costs}\label{cost}
We now consider the computational cost of the described methods, under the assumption that all the distances are available. Hence the cost for the construction of ${\mit \Gamma}$ is not accounted for.
When the graph is sparse,  adjacency lists are used for the graph representation.
 The cost to construct the similarity matrix $W$ is given by $\gamma=\gamma_1+\gamma_2+\gamma_3$, where

\smallskip\noindent
-- $\gamma_1$ is the cost of sparsifying ${\mit\Gamma}$,

\smallskip\noindent
-- $\gamma_2$ is the cost of computing the  scale parameters,

\smallskip\noindent
-- $\gamma_3$ is the cost of evaluating the similarity  function.

\medskip
For $F$-methods no sparsity is performed, hence $\gamma_1=0$ and $\gamma_3=O(n^2)$ since the similarity function  is evaluated on  a nearly full matrix.  For method $F_1$, $\gamma_2$ is of order $O(n^2)$ as it follows from Subsection \ref{mst}.  Also for methods  $F_2$ and $F_3$, $\gamma_2=O(n^2)$, indeed, given an integer $K < n$, $O(n)$ computational steps suffice to select the $K$th nearest neighbor to a point and to identify the set of the $h$th nearest neighbors, for $h=1,\ldots,  K$ (see \cite{cor}).

 For the same reason, the detection of the edges of ${\mit \Delta}$  for $E$-methods has a cost $\gamma_1=O(n^2)$.
 For method $E_1$ the costs  $\gamma_2$ and $\gamma_3$  are zero. The number $\eta$ of edges of ${\mit \Delta}$ is not a-priori quantifiable, and $\eta=O(n^2)$.  For method  $E_2$ $\gamma_2=O(\eta+n \log n)$  and  $\gamma_3=O(\eta)$. For methods $E_3$ and $E_4$ both $\gamma_2$ and $\gamma_3$ are of order $O(\eta)$. In any case all the costs result to be of order $O(n^2)$.

For $N$-methods and $M$-methods, the number of edges of ${\mit \Delta}$ is of order $O(K\,n)$ and their selection has a cost $\gamma_1$ of order $O(n^2)$.
When the unit similarity function is used, the costs  $\gamma_2$ and $\gamma_3$  are zero.
Otherwise, both   $\gamma_2$ and $\gamma_3$ are of order $O( K\,n)$.

\smallskip
In brief, all the above described methods have costs $\gamma$ of order $O(n^2)$. The methods which exploit  the sparsity of the graphs have computational costs dominated by the cost $\gamma_1$ of the sparsifying step.

\section{The experimentation}\label{expe}

The experiments have been conducted on an Intel(R)Core(TM)i7-4770CPU
@3.40GHz using double precision arithmetic.

\subsection{The datasets}\label{data}
Datasets belonging to four sets $A$, $B$, $R$, $U$ have been considered. The target clustering ${\mit \Pi}_{tar}$ is known for every dataset.

\begin {figure}[h!]
\begin{center}
\hskip 0.cm $A_1$ \hskip 2.8cm $A_2$\hskip 2.8cm $A_3$\hskip 2.8cm $A_4$\\[0.2cm]
\epsfig{file=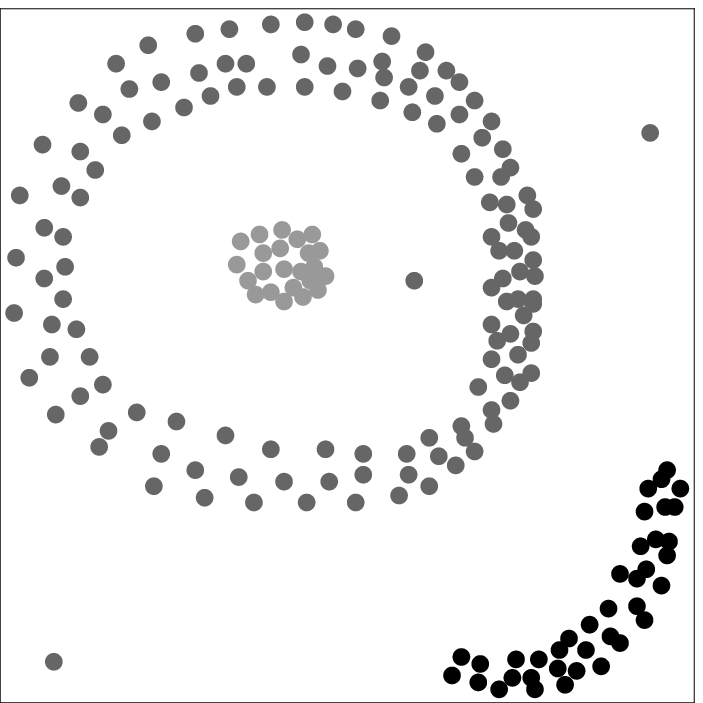,width=2.8cm} \quad
 \epsfig{file=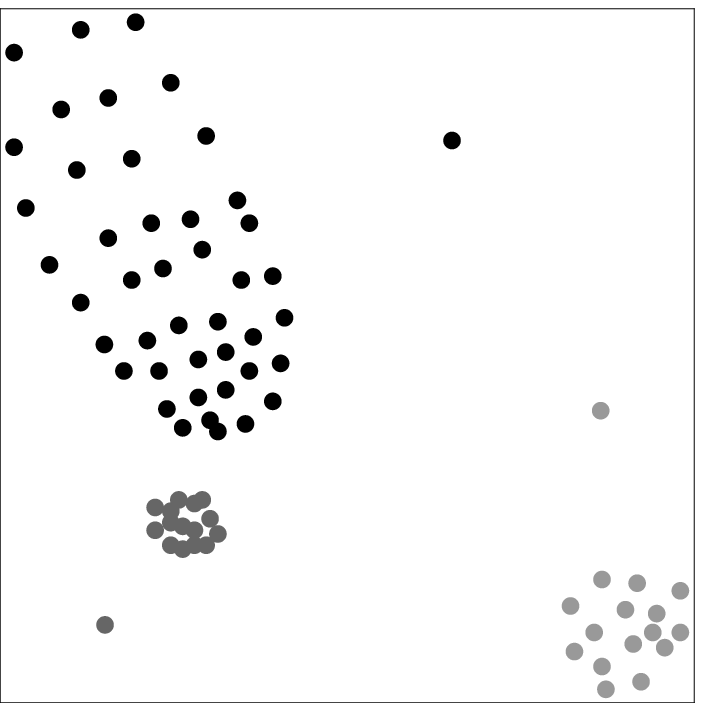,width=2.8cm}\quad
 \epsfig{file=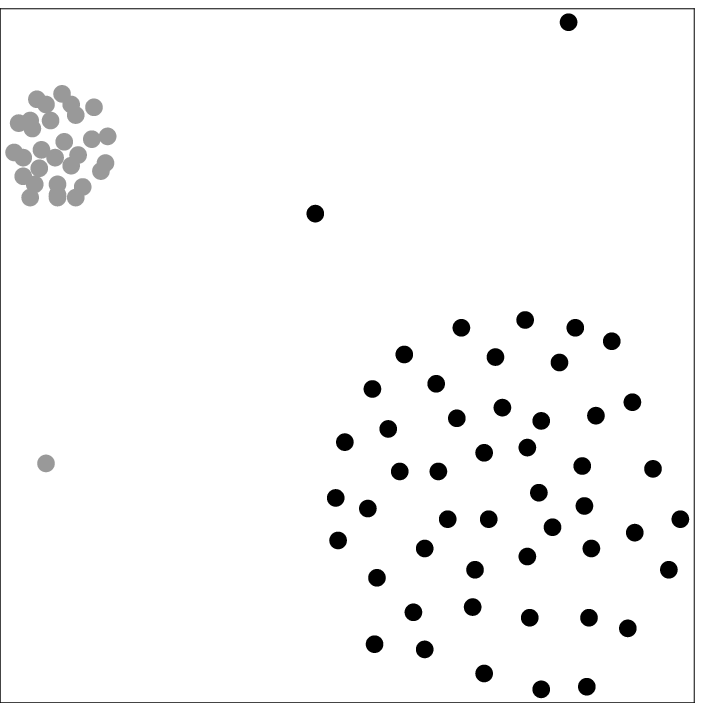,width=2.8cm} \quad
 \epsfig{file=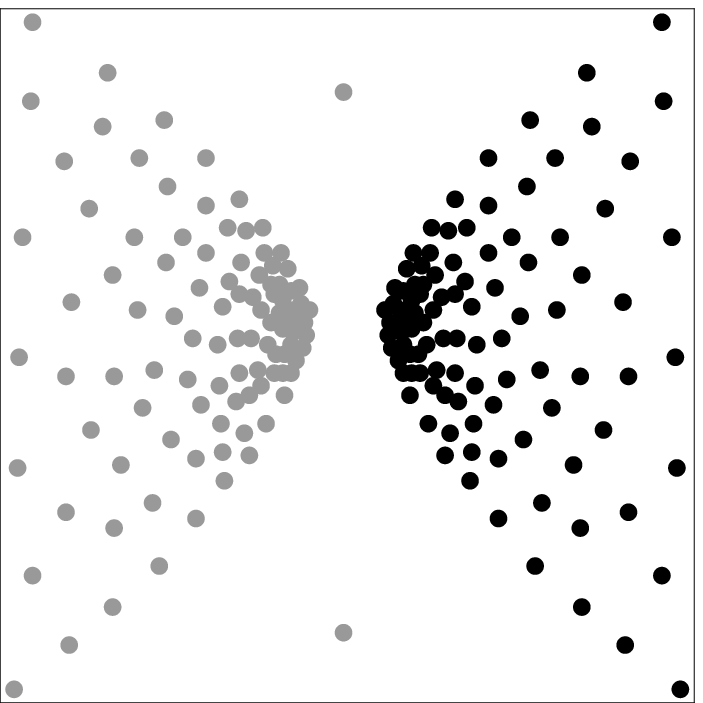,width=2.8cm}
 \end{center}
\begin{center}
\hskip 0.cm $A_5$ \hskip 2.8cm $A_6$\hskip 2.8cm $A_7$\hskip 2.8cm $A_8$\\[0.2cm]
\epsfig{file=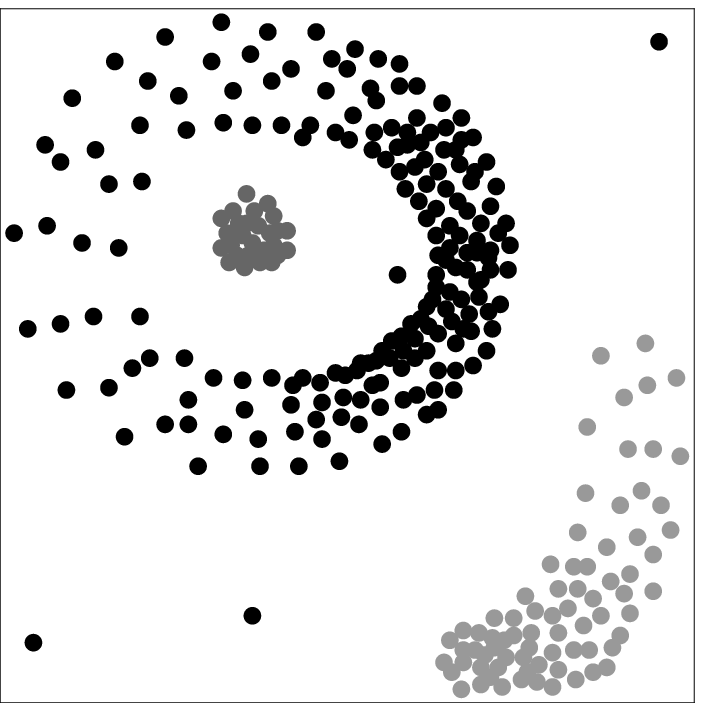,width=2.8cm} \quad
 \epsfig{file=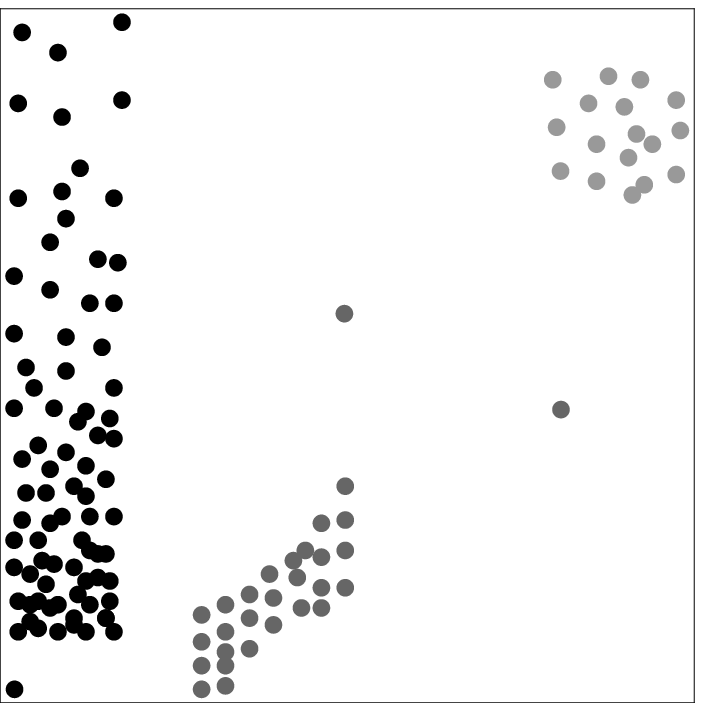,width=2.8cm}\quad
 \epsfig{file=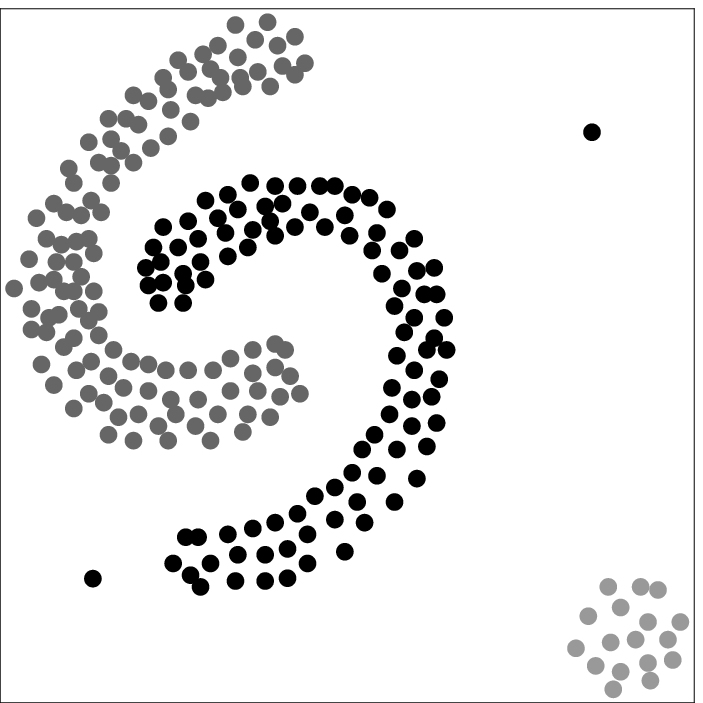,width=2.8cm} \quad
 \epsfig{file=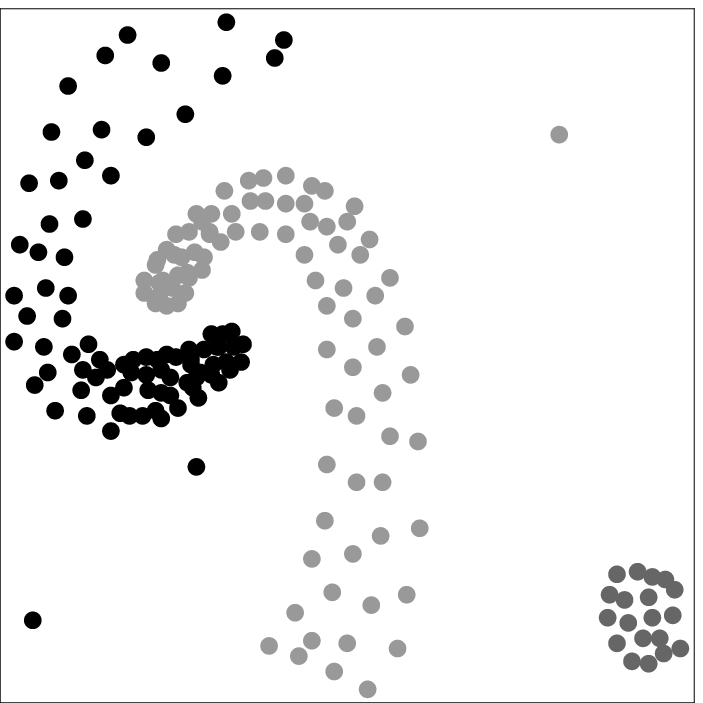,width=2.8cm}
 \end{center}
\begin{center}
\hskip 0.cm $A_9$ \hskip 2.7cm $B_1$\hskip 2.7cm $B_2$\hskip 2.7cm $B_3$\\[0.2cm]
\epsfig{file=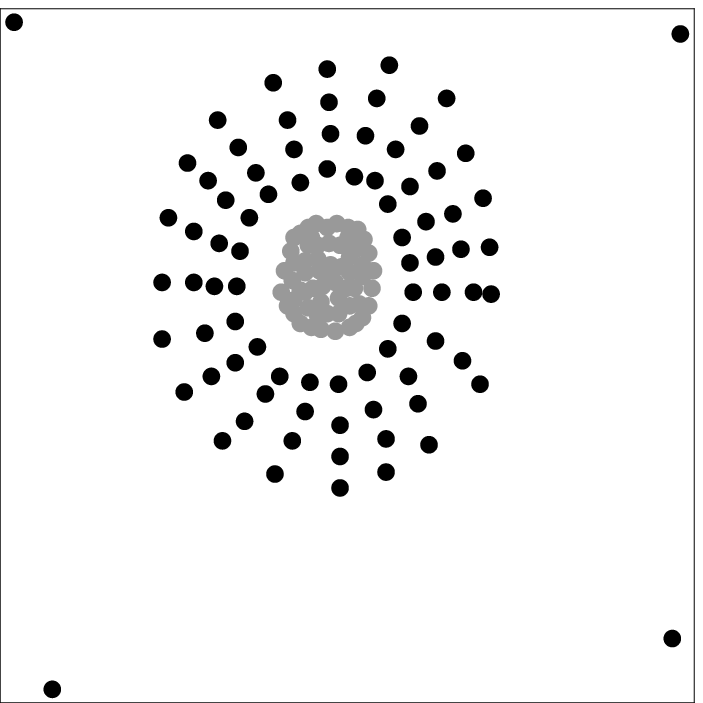,width=2.8cm} \quad
 \epsfig{file=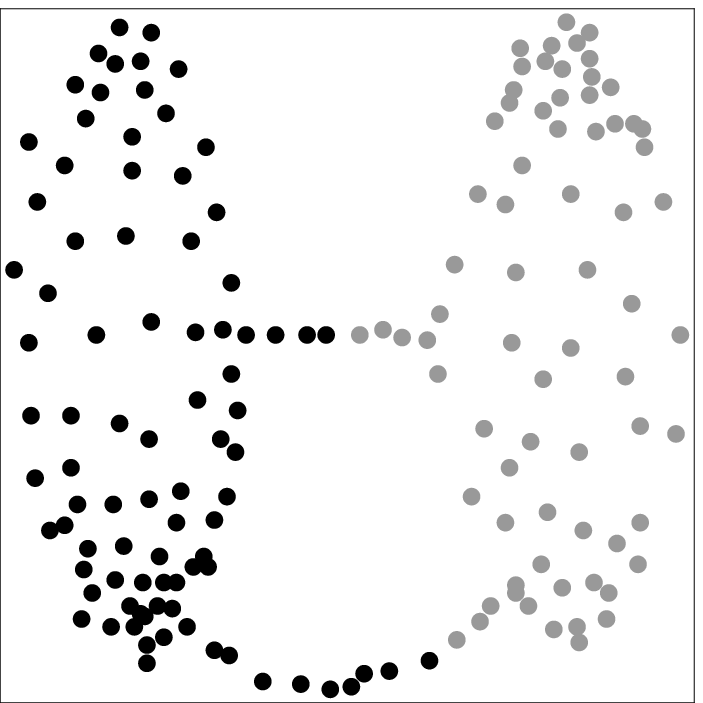,width=2.8cm}\quad
 \epsfig{file=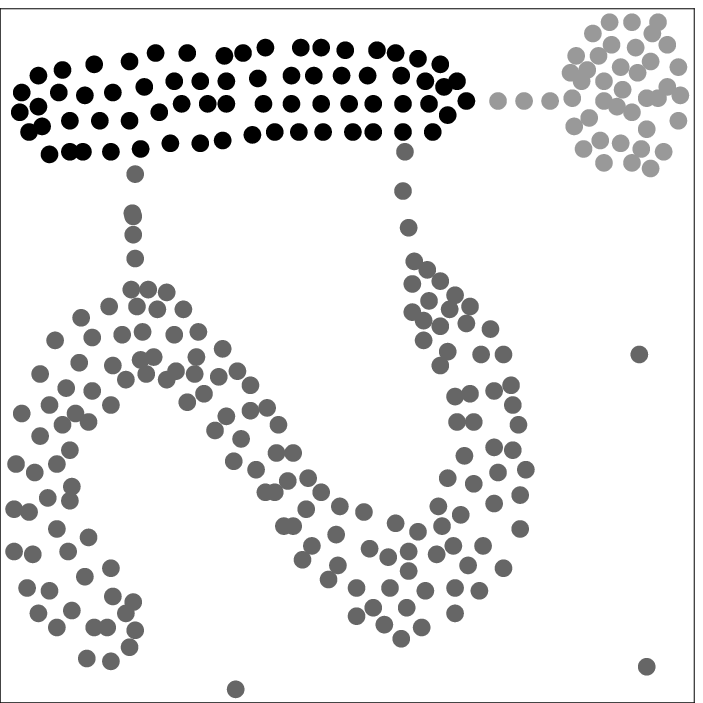,width=2.8cm} \quad
 \epsfig{file=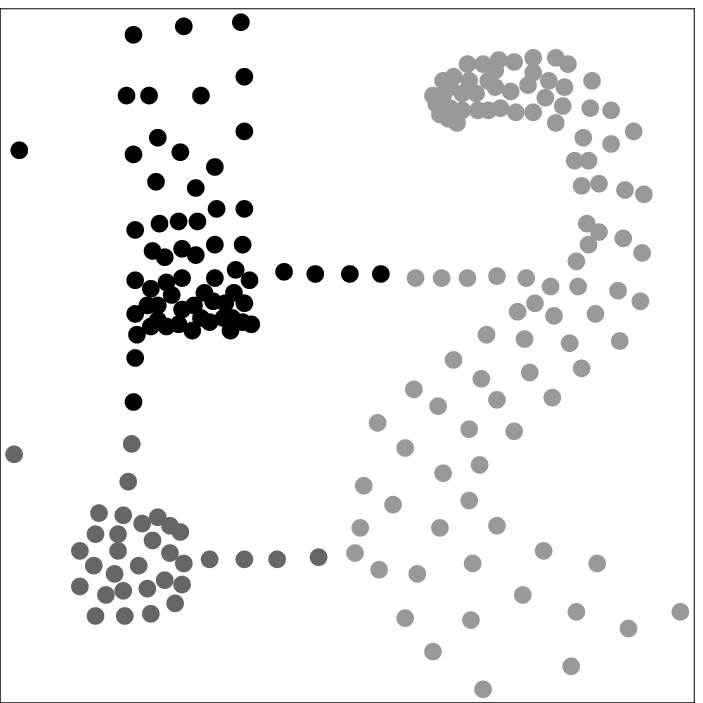,width=2.8cm}
 \end{center}
 \vskip -0.5cm
 \caption{\label{fig1} Artificial datasets generated in ${\bf R}^2$. }
\end{figure}

\smallskip\noindent
{\bf A}, {\bf B.}\quad Sets $A$ and $B$ contain 12 artificial datasets of ${\bf R}^2$ which can be found in the Spatial Data Mining Project of \cite{sou}. Their clusters have circular or noncircular shape, uniform or nonuniform density. In Figure \ref{fig1} they are shown according to growing difficulty. Different gray levels evidence the assumed target clusters. The first 9 datasets belong to set $A$ and have outliers and well separated clusters; the last 3 datasets belong to set $B$, have bridges which connect different clusters and lend themselves to possible different targets.
The datasets have from 76 to 289 points, normalized in such a way that
the maximum distance between points of each figure is equal to 1.

\begin {figure}[h!]
\begin{center}
$R_1$ \hskip 3.5cm $R_2$\hskip 3.5cm $R_3$\\
\epsfig{file=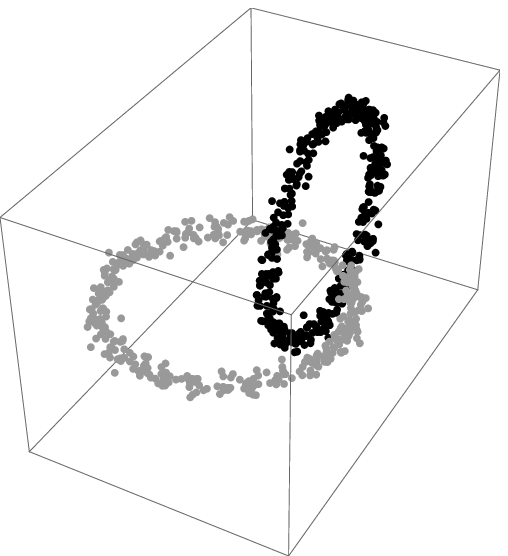,width=3.7cm}\quad
 \epsfig{file=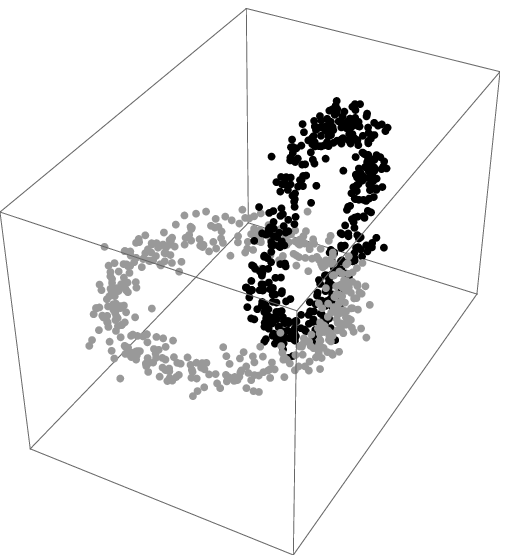,width=3.7cm}\quad
 \epsfig{file=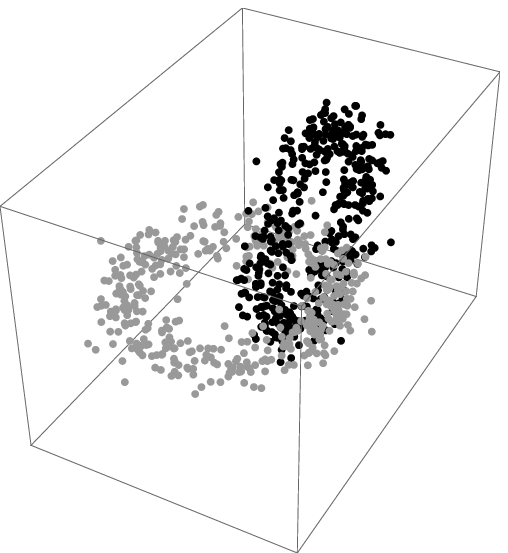,width=3.7cm}
 \end{center}
\vskip -0.5cm
 \caption{\label{fig2} Artificial datasets generated in ${\bf R}^3$. }
\end{figure}

\smallskip\noindent
{\bf R.}\quad Set $R$ contains 6 artificial datasets generated in ${\bf R}^3$ which consist of two interlacing rings with increasing data
dispersion and have 900 points each (see three of them in Figure
\ref{fig2}), and the 2 datasets generated in ${\bf R}^8$  with 200 points each which can be found in \cite{inka2} under the name I-$\Lambda$.

\smallskip\noindent
{\bf U.}\quad Set $U$ contains 4 real-world datasets {\it iris}, {\it wine}, {\it vote} and {\it seeds}
which have been taken from the UCI database \cite{lic}.

Dataset {\it iris}
is a typical test case for many statistical classification
techniques in machine learning. It consists of samples from three species of iris, with four features measured for each
sample: the length and the width of sepals and petals. Iris contains
$n=150$ points in ${\bf R}^4$ and 3 clusters are expected. Actually, the dataset contains only two
well separated clusters: one cluster contains the measurements of
only one species, while the other cluster contains the measurements
of the other two species and is not evidently separable without
further supervised information.

Dataset {\it wine} consists of samples of chemical analysis of
wines derived from three different cultivars. The analysis determined 13
attributes for each sample. Wine contains $n=178$ points in ${\bf R}^{13}$ and 3 clusters are expected.

Dataset {\it vote} consists of 435 samples. Each sample lists the
votes (a vote can be
"yes", "no" or can be missing) given by one U.S. congressman (out of 168 republican and 267
democrat congressmen) on 16 different questions. Vote contains $n=435$ points in ${\bf R}^{16}$ and 2 clusters are expected.

Dataset {\it seeds} consists of samples measuring the geometrical properties of the kernels of three different varieties of wheat: Kama, Rosa and Canadian, randomly selected for
the experiment. A soft X-ray technique was used to construct seven real-valued attributes. Seeds contains $n=210$ points in ${\bf R}^7$ and 3 clusters are expected.

\subsection{Results}\label{res}
The $q$ methods listed in Subsection \ref{def} have been applied to each dataset obtaining almost always clusterings ${\mit \Pi}_{obt}$ with the expected number of clusters.
For $r=1,\ldots,q$ and $s=1,\ldots,\mu$, the $r$th method has been applied to the $s$th dataset, with $\mu=9$ for the set $A$, $\mu=3$ for the set $B$, $\mu=8$ for the set $R$, $\mu=4$ for the set $U$. The following quantities have been computed:\\
-- the sparsity level $\theta_{r,s}$ of the similarity graph $G$,\\
-- the accuracy $\alpha_{r,s}$, measured by the NMI index of the clustering ${\mit \Pi}_{obt}$.

\noindent
For each set $A$, $B$, $R$ and $U$, the averaged sparsity level and the averaged accuracy of the $r$th method, $r=1,\ldots,q,$
\[
\widetilde \theta_r=\frac{1}{\mu}\, \sum_{s=1}^\mu \theta_{r,s}
\quad\hbox{and}\quad
\widetilde\alpha_r=\frac{1}{\mu}\, \sum_{s=1}^\mu \alpha_{r,s}
   \]
are computed varying the $\mu$ datasets in the set.

The sparsity level is measured with respect to the machine zero $eps= 2^{-52}$, i.e.
$\theta_{r,s}=\nu/n^2$, where
$\nu$ is the number of elements in $W$ which are smaller than $eps$. When the Gaussian similarity  function is used, the reference to $eps$ instead of the real zero makes significant sparsity levels even in the case of a (theoretical) dense graph, as can be seen in Table \ref{ta1}. Except for the $F$-methods, $\widetilde\theta_r$ vary negligibly when the methods belong to the same class, pointing out that once a sparsity model has been chosen, the number of edges in ${\mit \Delta}$ is almost unaffected by the choice of $\sigma$. The index $i$ stands for any $i=1,\dots,4$.

\begin{table} [!h]
\centering
\renewcommand{\arraystretch}{0.4}
\scriptsize
\begin{tabular}{|c|c|cccc|}
\hline
&&&&&\\
&method&{Set A}&{Set B}&{Set R}&{Set U}\\[0.15cm]
\hline
&&&&& \\
&$F_1$ & 52\% &  74\% &  67\% &  44\%  \\[0.2cm]
&$F_2$ & 19\% &  24\% &  40\% &  24\%  \\[0.2cm]
&$F_3$ & 16\% &  18\% &  39\% &  18\% \\[0.2cm]
\hline
&&&&&   \\
&$E_i$ &  92\% &  95\% & 93\%  & 94\%   \\[0.2cm]
$K=1+\lfloor\log_2 n\rfloor$&$N_i$ &  94\% &  95\%  &  98\% &  95\%  \\[0.2cm]
&$M_i$ &  96\% &  97\%  &  99\%  &  97\%  \\[0.2cm]
\hline
&&&&&   \\
&$E_i$ &  88\%  &  92\%  & 90\%  & 91\%   \\[0.2cm]
$K=1+\lfloor\sqrt n\rfloor$&$N_i$ &  90\% &  91\%  &  94\%  &  91\%  \\[0.2cm]
&$M_i$ &  93\% &  94\%  &  97\%  &  95\%  \\[0.2cm]
\hline
\end{tabular}
\caption{ \label{ta1} Averaged sparsity levels of the methods.}
\end{table}

Let $\rho_r$ denote the rank of the $r$th method, i.e. the position of $\widetilde\alpha_r$, in the not increasing sequence $\widetilde{\bm \alpha}=\{\widetilde\alpha_1,\ldots,\widetilde\alpha_q\}$. Methods with equal accuracy receive the same ranking number according to the standard competition ``1224'' policy (i.e. the ranking of each method is equal to 1 plus the number of the methods which have a larger accuracy).
Table \ref{ta2} shows the averaged accuracies $\widetilde \alpha_r$ and the corresponding ranking $\rho_r$.
From Tables  \ref{ta1} and  \ref{ta2}, it appears that there is not a strict relation between the sparsity level of a method and its accuracy, in the sense that the accuracy seems to depend on ``which'' more than on ``how many'' edges are retained in the graph ${\mit \Delta}$. For example on class $A$, the $M$-methods with both values of $K$ get good results while the $N$-methods with $K = K_\ell$ which  have an intermediate sparsity level, obtain poor results.

\begin {figure}[h!]
\begin{center}
\epsfig{file=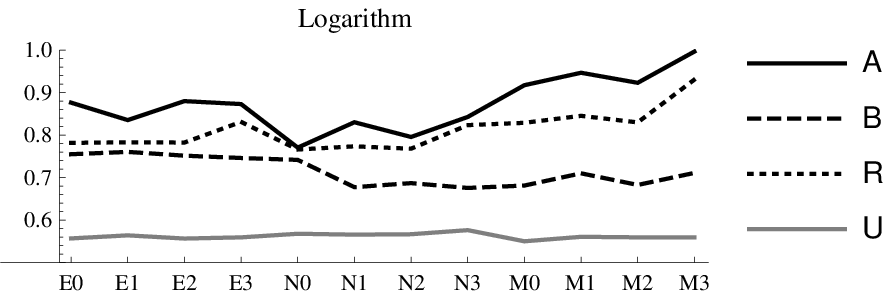,width=10cm}
 \end{center}
 \caption{\label{graf1} Averaged accuracies for $E$, $N$ and $M$-methods. }
\end{figure}

\begin {figure}[h!]
\begin{center}
\epsfig{file=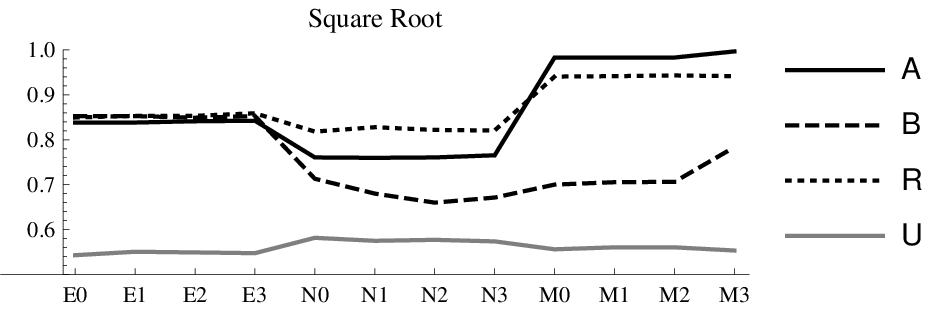,width=10cm}
 \end{center}
 \caption{\label{graf2} Averaged accuracies for $E$, $N$ and $M$-methods. }
\end{figure}

\begin{table} [!h]
\centering
\renewcommand{\arraystretch}{0.3}
\scriptsize
\begin{tabular}{|c|c|cc|cc|cc|cc|}
\hline
&&&&&&&&&\\
&&\multicolumn{2}{|c|}{Set A}&\multicolumn{2}{|c|}{Set B}&\multicolumn{2}{|c|}{Set R}&\multicolumn{2}{|c|}{Set U}\\
&&&&&&&&& \\
\cline{3-10}
&&&&&&&&& \\
& method&$\widetilde \alpha$& $\rho$
 &$\widetilde \alpha$& $\rho$
 &$\widetilde \alpha$& $\rho$
 &$\widetilde \alpha$& $\rho$\\
&&&&&&&&& \\
\hline
\hline
&&&&&&&&& \\
&$F_1$ & 0.87 & 11 &  0.59 & 27 &  0.84 & 11 &  0.52 & 27  \\[0.15cm]
&$F_2$ & 0.73 & 26 &  0.69 & 18 &  0.76 & 26 &  0.58 & 1  \\[0.15cm]
&$F_3$ & 0.72 & 27 &  0.72 & 11 & 0.72 & 27 &  0.58 & 1  \\[0.15cm]
\hline
\hline
&&&&&&&&&   \\
&$E_1$ & 0.88 & 9 & 0.76& 6 & 0.78 & 20 &  0.56 & 11  \\[0.15cm]
&$E_2$ &  0.84 & 13 &  0.76 & 6 & 0.78 & 20 & 0.56 &  11  \\[0.15cm]
&$E_3$ &  0.88 &9 &  0.75 & 8 &  0.78 & 20 & 0.56 & 11  \\[0.15cm]
&$E_4$ &  0.87 &11 &  0.75 & 8 & 0.83 & 12 & 0.56 & 11  \\[0.15cm]
\cline{2-10}
&&&&&&&&& \\
&$N_1$ & 0.77&21 &  0.74 & 10 &  0.77 & 23 &  0.57 & 7  \\[0.15cm]
$K=1+\lfloor\log_2 n\rfloor$&$N_2$ &  0.83&19 &  0.68 & 20 &  0.77 & 23 &  0.57 & 7 \\[0.15cm]
&$N_3$ &  0.80&20 &  0.69 & 18 &   0.77 & 23 &  0.57 & 7  \\[0.15cm]
&$N_4$ &  0.84&13 &  0.68 & 20 &  0.82 & 16 &  0.58 & 1 \\[0.15cm]
\cline{2-10}
&&&&&&&&& \\
&$M_1$ &  0.92&7 &  0.68 & 20 & 0.83 &  12 &  0.55 & 21 \\[0.15cm]
&$M_2$ &  0.95&6 &  0.71 & 12 &  0.85 & 7 &  0.56 & 11 \\[0.15cm]
&$M_3$ &  0.92&7 &  0.68 & 20 & 0.83 & 12 &   0.56 & 11 \\[0.15cm]
&$M_4$ &  1.00&1 &  0.71 & 12 &  0.93 & 5 &  0.56 & 11 \\[0.15cm]
\hline
\hline
&&&&&&&&&   \\
&$E_1$ & 0.84 & 13 & 0.85& 1 & 0.85 &7 &  0.54 & 26  \\[0.15cm]
&$E_2$ &  0.84 & 13 &  0.85 & 1 & 0.85 & 7 & 0.55 &  21  \\[0.15cm]
&$E_3$ &  0.84 &13 &  0.85 & 1 &  0.85 & 7 & 0.55 & 21  \\[0.15cm]
&$E_4$ &  0.84 &13 &  0.85 & 1 & 0.86 & 6 & 0.55 & 21  \\[0.15cm]
\cline{2-10}
&&&&&&&&& \\
&$N_1$ & 0.76& 23 &  0.71 & 12 &  0.82 & 16 &  0.58 & 1 \\[0.15cm]
$K=1+\lfloor\sqrt n\rfloor$&$N_2$ &  0.76&23 &  0.68 & 20 &  0.83 & 12 &  0.58 & 1 \\[0.15cm]
&$N_3$ &  0.76&23 &  0.66 & 26 &  0.82 & 16 &   0.58 & 1 \\[0.15cm]
&$N_4$ &  0.77&21 &  0.67 & 25 &  0.82 & 16 &   0.57 & 7 \\[0.15cm]
\cline{2-10}
&&&&&&&&& \\
&$M_1$ & 0.98&3 &  0.70 & 17 & 0.94 & 1 &  0.56 & 11 \\[0.15cm]
&$M_2$ &  0.98&3 &  0.71 & 12 &  0.94 & 1 &  0.56 & 11 \\[0.15cm]
&$M_3$ &  0.98&3 &  0.71 & 12 & 0.94 & 1 &  0.56 & 11 \\[0.15cm]
&$M_4$ &  1.00&1 &  0.78 & 5 &  0.94 & 1 &  0.55 & 21 \\[0.15cm]
\hline

\end{tabular}
\caption{ \label{ta2} Averaged accuracies $\widetilde \alpha$ and ranking of the methods.}
\end{table}
Figures \ref{graf1} and \ref{graf2} plot the averaged accuracies shown in Table \ref{ta2}, for $E$, $N$ and $M$-methods, grouped according to the value of $K$. From these figures it appears evident that the choice of the similarity function   influences more the averaged accuracy, when $K =K_\ell$.


It is worth noting that the averaged accuracy of all the different methods for problems of class $U$ is very low and the range of its  values is very small. This class of problems is not worthwhile to analyze and compare the considered methods. Restricting ourselves to the classes of problems $A$, $B$ and $R$, from Table \ref{ta2} it appears that the $F$-methods and  the $N$-methods are outperformed by $M$-methods on problems $A$ and $R$ and  by $E$-methods on problems $B$. For this reason we limit our analysis to $E$ and $M$-methods.

 Since the averaged values of the accuracy are not sufficient to show the different features of these methods, we analyze the behavior of selected methods on selected problems. We choose a representative method for both $E$ and $M$-methods. More precisely, due to their good performance on average,  the methods $E_4$ and $M_4$, for the two values  $K_\ell$ and $K_s$ are selected.

\begin{table} [!h]
\centering
\renewcommand{\arraystretch}{0.4}
\scriptsize
\begin{tabular}{|c|c|cccccccc|}
\hline
&&&&&&&&&\\
&method&$ A_7$&$ A_8$&$ A_9$&$ B_1$&$ B_2$&$ B_3$&$ (I-\Lambda)_1$&$ (I-\Lambda)_2$\\[0.15cm]
\hline
&&&&&&&&&\\
$K=1+\lfloor\log_2 n\rfloor$&$E_4$ &  1 &  0.54  &  0.35 &  0.38 &  0.9 &  0.96 &  0.52  &  0.24 \\[0.2cm]
&$M_4$ &  1 &  1  &  1 &  0.91 &  0.47 &  0.75 &  0.81  &  0.81\\[0.2cm]
\hline
&&&&&&&&&\\
$K=1+\lfloor\sqrt n\rfloor$&$E_4$ &  0.68 &  0.58  &  0.35 &  0.85 &  0.9 &  0.8 &  0.57  &  0.51  \\[0.2cm]
&$M_4$ &  1 &  1  &  1 &  0.95 &  0.47 &  0.93 &  0.81  &  0.81 \\[0.2cm]
\hline
\end{tabular}
\caption{ \label{ta3} Accuracies of methods $E_4$ and $M_4$ on selected problems.}
\end{table}

\begin {figure}[h!]
\begin{center}
\epsfig{file=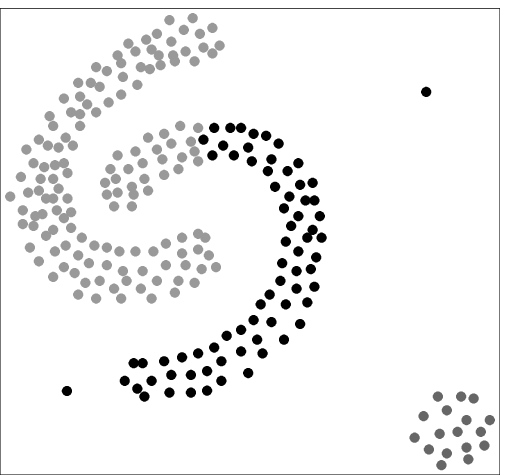,width=3.5cm}\quad
 \epsfig{file=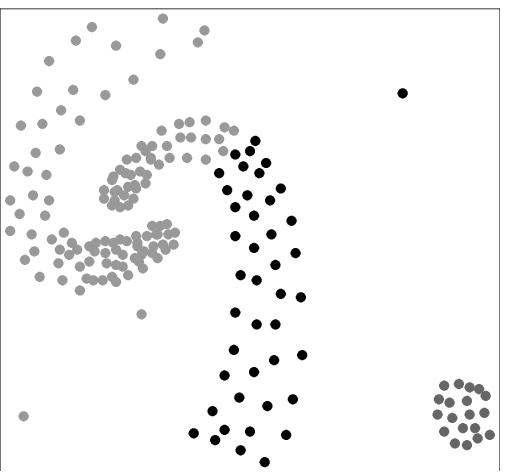,width=3.5cm}\quad
 \epsfig{file=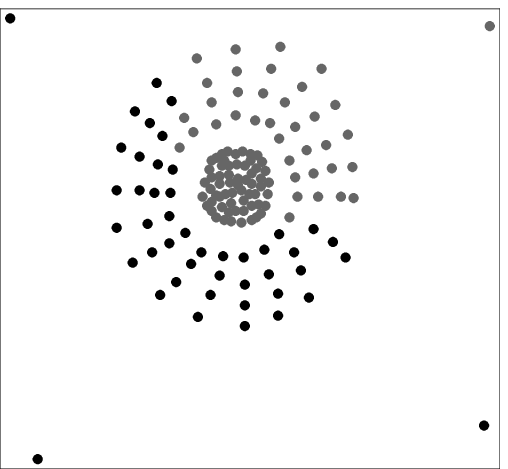,width=3.5cm} \quad\\[0.2cm]
 \hskip 0.cm $\alpha=0.68$ \hskip 2.5cm $\alpha=0.58$\hskip 2.5cm $\alpha=0.35$
 \end{center}
 \caption{\label{ponti0}  Clustering obtained for problems $A_7$ (left), $A_8$ (middle) and $A_9$ (right) by applying method $E_4$, with $K = K_s$.}
 \end{figure}

\begin {figure}[h!]
\begin{center}
 \epsfig{file=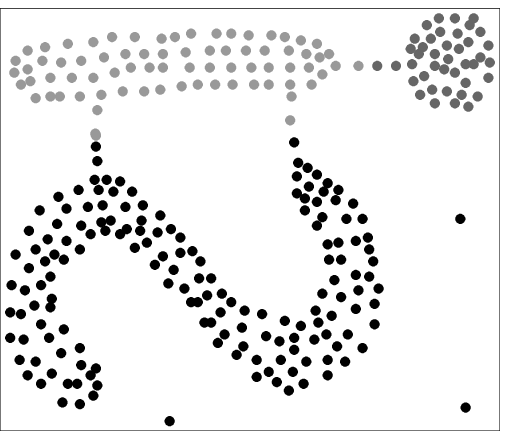,width=3.5cm}\quad
 \epsfig{file=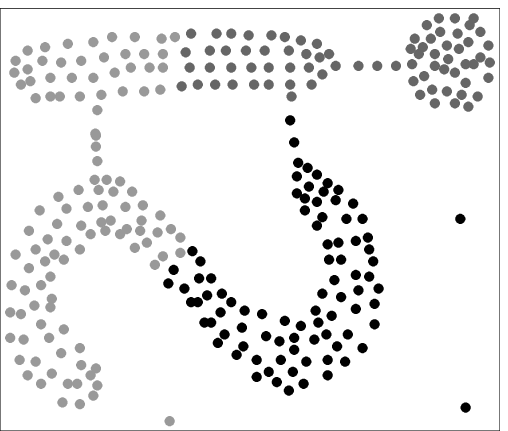,width=3.5cm} \quad\\[0.2cm]
 \hskip 0.cm $\alpha=0.9$ \hskip 2.5cm $\alpha=0.47$
 \end{center}
 \caption{\label{ponti} Clustering obtained for problem $B_2$ by applying method $E_4$ (left),  and $M_4$ (right)  with $K = K_\ell$.}
 \end{figure}

\begin {figure}[h!]
\begin{center}
 \epsfig{file=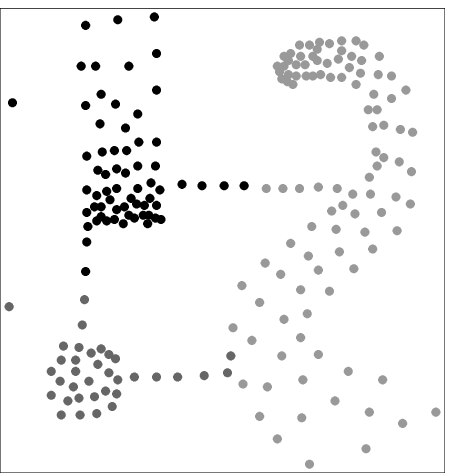,width=3.5cm}\quad
 \epsfig{file=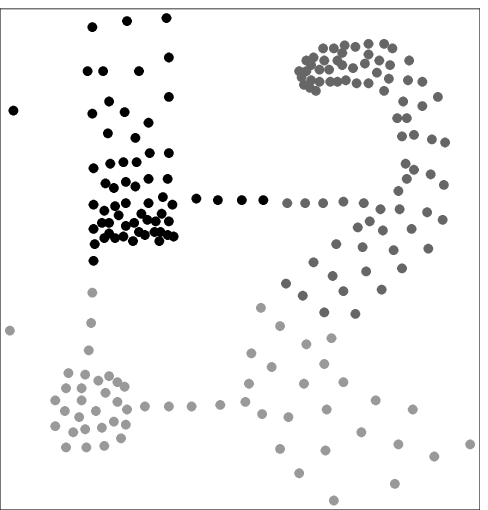,width=3.5cm} \quad\\[0.2cm]
 \hskip 0.cm $\alpha=0.96$ \hskip 2.5cm $\alpha=0.75$
 \end{center}
 \caption{\label{ponti1}  Clustering obtained for problem $B_3$ by applying method $E_4$ (left),  and $M_4$ (right)  with $K = K_\ell$.}
 \end{figure}

\begin {figure}[h!]
\begin{center}
 \epsfig{file=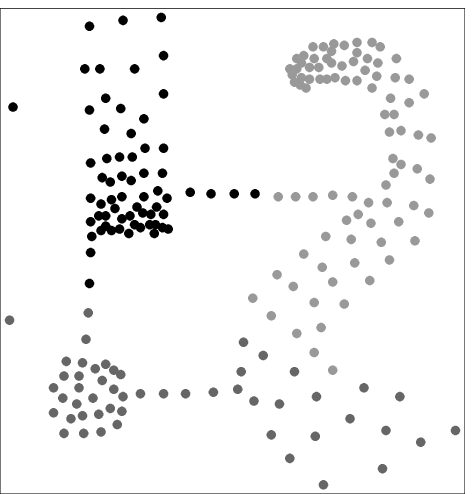,width=3.5cm}\quad
 \epsfig{file=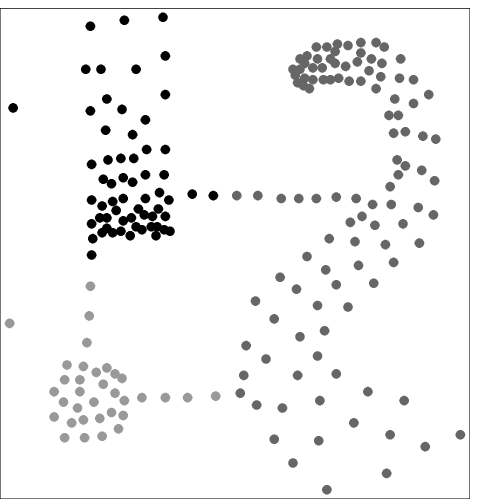,width=3.5cm} \quad\\[0.2cm]
 \hskip 0.cm $\alpha=0.8$ \hskip 2.5cm $\alpha=0.93$
 \end{center}
 \caption{\label{ponti2}  Clustering obtained for problem $B_3$ by applying method $E_4$ (left),  and $M_4$ (right)  with $K = K_s$.}
 \end{figure}

 In Table \ref{ta3} we list the accuracies of the four methods under investigation on three problems of class $A$, namely $A_7$, $A_8$, $A_9$, on all problems of class $B$ and on the two $I-\Lambda$ problems of class $R$. On all the remaining problems of class $A$ and $R$ these four methods find the target clusters. By inspection of Table \ref{ta3}, it appears that the $M$-methods always outperform $E$-methods except for problems with bridges. For $E$-methods the choice of the radius of the neighbors can be crucial especially when there are clusters with different densities. In order to better understand the meaning of Table \ref{ta3}, see Figures \ref{ponti0}, \ref{ponti}, \ref{ponti1} and \ref{ponti2} where  critical clustering outcomes are plotted and the corresponding values of accuracy $\alpha$ are given.

Figure \ref{ponti0} shows the poor performance of  method $E_4$ with $K = K_s$ on problems $A_7$, $A_8$ and $A_9$. While for problem $A_7$ the choice of a smaller radius allows finding the target clustering, for problems $A_8$ and $A_9$  also the choice  $K = K_\ell$ produces analogous results, may be due to nonuniform densities.
 Figure  \ref{ponti} shows the clustering obtained by   method $E_4$ (on the left) having irrelevant differences from the target one and the wrong clustering obtained by  method $M_4$ (on the right) for problem $B_2$ which has short bridges. Similar results are obtained for $K = K_s$.
Figures \ref{ponti1} and \ref{ponti2} show that the choice of $K$ is critical for both methods $E_4$ and $M_4$   on problem  $B_3$ which has bridges and clusters of nonuniform densities.
 To further investigate the negative influence of the bridges, we have removed bridges from the problems of class $B$, obtaining the   target clusterings with all the four methods under consideration. In this case the nonuniform densities of problems $B_1$ and $B_3$ do not affect the behavior of $E$-methods, due to the sufficiently large distances among the clusters.

\subsection{Conclusions}
In this paper several techniques to construct the similarity matrix have been considered. The corresponding methods have been tested on different datasets by applying a normalized spectral algorithm. The methods based on the dense graph and those obtained by exploiting the non mutual $K$-nearest neighbor sparsity are clearly outperformed by those based on the $\epsilon$-neighbor sparsity and mutual $K$-nearest neighbor sparsity. The performance of the more effective $E$ and $M$-methods depends on the characteristics of the problems. Among the different $M$-methods, the one with  $K = K_s$ and Gaussian similarity function given in (\ref{gau}) with $\sigma = s$ given in (\ref{sigma}) appears to be the most reliable, even if in few cases it is outperformed by the corresponding $E$-method.

\end{document}